\documentclass[11pt]{amsart} 
\usepackage{amsmath,amssymb,amsthm,mathrsfs}
\newtheorem{theorem}{Theorem}
\newtheorem{proposition}[theorem]{Proposition}
\newtheorem{corollary}[theorem]{Corollary}
\newtheorem{lemma}[theorem]{Lemma}

\begin{document}

\title{Embedded Weingarten tori in $S^3$} 
\author{Simon Brendle}
\begin{abstract}
In this paper, we show that an embedded Weingarten surface in $S^3$ of genus $1$ must be rotationally symmetric, provided that certain structure conditions are satisfied. The argument involves an adaptation of our proof of Lawson's Conjecture for minimal tori.
\end{abstract}
\address{Department of Mathematics \\ Stanford University \\ Stanford, CA 94305}
\thanks{The author was supported in part by the National Science Foundation under grant DMS-1201924.}
\maketitle

\section{Introduction}

In a recent paper \cite{Brendle1}, we proved that the Clifford torus is the only embedded minimal surface in $S^3$ of genus $1$, thereby giving an affirmative answer to a question posed by Lawson \cite{Lawson2} in 1970 (see also \cite{Brendle-survey}). The argument in \cite{Brendle1} can be extended in various ways. For example, in \cite{Brendle2}, we showed that any minimal torus which is immersed in the sense of Alexandrov must be rotationally symmetric. Moreover, it was observed in \cite{Andrews-Li} that the arguments in \cite{Brendle1} can be extended to the setting of constant mean curvature tori in $S^3$. 

In this note, we extend the arguments in \cite{Brendle1} to a class of Weingarten tori in $S^3$. Let $F: \Sigma \to S^3$ be a surface in $S^3$ which satisfies a PDE of the form 
\[\lambda_1+\lambda_2 = \psi(\lambda_1-\lambda_2),\] 
where $\lambda_1 \geq \lambda_2$ denote the principal curvatures. A surface with this property is referred to as a Weingarten surface. 

\begin{theorem}
\label{main.result}
Let $\psi(s)$ be an even function satisfying $0 \leq s \, \psi'(s) < \min \{\psi(s),s\}$ and $0 \leq s \, \psi''(s) \leq 1-\psi'(s)^2$ for $s \geq 0$. Let $F: \Sigma \to S^3$ be an embedded torus which satisfies the equation $\lambda_1 + \lambda_2 = \psi(\lambda_1-\lambda_2)$. Then $F$ is rotationally symmetric. More precisely, we can find an anti-symmetric matrix $Q \in \mathfrak{so}(4)$ of rank $2$ such that $Q \, F(x) \in \text{\rm span} \{\frac{\partial F}{\partial x_1}(x),\frac{\partial F}{\partial x_2}(x)\}$ for all $x \in \Sigma$.
\end{theorem}

For each $a \in (0,\infty)$, each $b \in [0,1]$, and each $c \in [0,\infty)$, the function $\psi(s) = \sqrt{a + b \, s^2} + c$ satisfies the assumptions of Theorem \ref{main.result}. Hence, we can draw the following conclusion: 

\begin{corollary}
Let us fix real numbers $a \in (0,\infty)$, $b \in [0,1]$, and $c \in [0,\infty)$. Moreover, let $F: \Sigma \to S^3$ be an embedded torus satisfying $(\lambda_1+\lambda_2-c)^2 = a + b \, (\lambda_1-\lambda_2)^2$ and $\lambda_1+\lambda_2 \geq c$. Then $F$ is rotationally symmetric in the sense described above. 
\end{corollary}

The proof of Theorem \ref{main.result} is similar in spirit to our earlier work \cite{Brendle1} on the Lawson conjecture, and involves an application of the maximum principle to a two-point function. This technique was pioneered by Huisken in \cite{Huisken}, and was developed further by Andrews \cite{Andrews} (see also \cite{Andrews-Langford-McCoy}).

Instead of assuming that $F$ is an embedding, it suffices to assume that $F$ is an Alexandrov immersion. This can be seen by adapting the arguments in \cite{Brendle2}. Finally, we note that several rigidity results for Weingarten surfaces have been established in \cite{Brendle-Eichmair}, \cite{Bryant}, \cite{Hopf}, \cite{Montiel-Ros}.

\section{Absence of umbilical points}

In this section, we show that a Weingarten torus satisfying the assumptions of Theorem \ref{main.result} is free of umbilical points. This is motivated by the classical result of Lawson \cite{Lawson1} that a minimal torus in $S^3$ has no umbilical points. The proof is an adaptation of ideas in Bryant's paper \cite{Bryant}.

\begin{proposition}
\label{absence.of.umbilic.points}
Let $F: \Sigma \to S^3$ be an immersed torus which satisfies the equation $\lambda_1 + \lambda_2 = \psi(\lambda_1-\lambda_2)$. Then $F$ has no umbilical points.
\end{proposition}

\textbf{Proof.} 
By the uniformization theorem, we can assume that $F$ is a conformal immersion from a flat torus $\Sigma = \mathbb{C} / \Lambda$ into $S^3$. Let $z$ be the standard complex coordinate on $\mathbb{C}$. Since $F$ is conformal, we may write $g_{z\bar{z}} = 2 \, e^{2\rho}$ for some smooth function $\rho$. It is straightforward to verify that  
\[h_{z\bar{z}} = e^{2\rho} \, (\lambda_1+\lambda_2)\] 
and 
\[|h_{zz}| = e^{2\rho} \, (\lambda_1-\lambda_2).\] 
Therefore, the equation $\lambda_1 + \lambda_2 = \psi(\lambda_1-\lambda_2)$ can be rewritten as 
\[e^{-2\rho} \, h_{z\bar{z}} = \psi(e^{-2\rho} \, |h_{zz}|).\] 
We now differentiate this identity with respect to $z$. Since $\psi(s)$ is an even function, we may write $\psi'(s) = s \, \chi(s)$ for some Lipschitz continuous function $\chi(s)$. Using the Codazzi equations, we obtain 
\begin{align*} 
e^{-2\rho} \, \partial_{\bar{z}} h_{zz} 
&= \partial_z (e^{-2\rho} \, h_{z\bar{z}}) \\ 
&= \psi'(e^{-2\rho} \, |h_{zz}|) \, \partial_z(e^{-2\rho} \, |h_{zz}|) \\ 
&= \frac{1}{2} \, \chi(e^{-2\rho} \, |h_{zz}|) \, \partial_z(e^{-4\rho} \, |h_{zz}|^2) \\ 
&= \frac{1}{2} \, \chi(e^{-2\rho} \, |h_{zz}|) \, e^{-4\rho} \, (h_{\bar{z}\bar{z}} \, \partial_z h_{zz} + \overline{h_{\bar{z}\bar{z}} \, \partial_{\bar{z}} h_{zz}} - 4 \, \partial_z \rho \, |h_{zz}|^2). 
\end{align*} 
Therefore, the function $h_{zz}$ satisfies the Beltrami-type equation 
\[\partial_{\bar{z}} h_{zz} = \frac{1}{2} \, \chi(e^{-2\rho} \, |h_{zz}|) \, e^{-2\rho} \, (h_{\bar{z}\bar{z}} \, \partial_z h_{zz} + \overline{h_{\bar{z}\bar{z}} \, \partial_{\bar{z}} h_{zz}} - 4 \, \partial_z \rho \, |h_{zz}|^2).\] 
In particular, the zeroes of the function $h_{zz}$ are isolated and have finite order (see \cite{Bers-John-Schechter}, p.~259). Moreover, the Hopf indices at the zeroes of the function $h_{zz}$ all have the same sign. Since $\Sigma$ is a torus, a standard index count argument implies that $h_{zz}$ cannot have any zeroes at all. \\

\section{The main calculation}

Let $F: \Sigma \to S^3$ be an embedding, and let $\nu$ denote a choice of unit normal vector field along $F$. In other words, for each point $x \in \Sigma$, the vector $\nu(x)$ is normal to the surface, but tangential to $S^3$. We assume that the principal curvatures of $F$ satisfy the equation 
\[\lambda_1+\lambda_2 = \psi(\lambda_1-\lambda_2),\] 
where $\psi$ satisfies the structure conditions listed in Theorem \ref{main.result}. 

We next consider a positive function $\Phi$ on $\Sigma$ with the property that 
\[Z(x,y) = \Phi(x) \, (1 - \langle F(x),F(y) \rangle) + \langle \nu(x),F(y) \rangle \geq 0\] 
for all points $x,y \in \Sigma$. Finally, we assume that there exists a pair of points $\bar{x} \neq \bar{y}$ satisfying $Z(\bar{x},\bar{y}) = 0$. Let $(x_1,x_2)$ be geodesic normal coordinates around $\bar{x}$, and let $(y_1,y_2)$ be geodesic normal coordinates around $\bar{y}$. 

At the point $(\bar{x},\bar{y})$, we have 
\begin{align*} 
0 = \frac{\partial Z}{\partial x_i}(\bar{x},\bar{y}) 
&= \frac{\partial \Phi}{\partial x_i}(\bar{x}) \, (1 - \langle F(\bar{x}),F(\bar{y}) \rangle) \\ 
&- \Phi(\bar{x}) \, \Big \langle \frac{\partial F}{\partial x_i}(\bar{x}),F(\bar{y}) \Big \rangle + h_i^k(\bar{x}) \, \Big \langle \frac{\partial F}{\partial x_k}(\bar{x}),F(\bar{y}) \Big \rangle 
\end{align*}
and 
\[0 = \frac{\partial Z}{\partial y_i}(\bar{x},\bar{y}) = -\Phi(\bar{x}) \, \Big \langle F(\bar{x}),\frac{\partial F}{\partial y_i}(\bar{y}) \Big \rangle + \Big \langle \nu(\bar{x}),\frac{\partial F}{\partial y_i}(\bar{y}) \Big \rangle.\] 
Without loss of generality, we may assume that the second fundamental form at $\bar{x}$ is diagonal, so that $h_{11}(\bar{x}) = \lambda_1$, $h_{12}(\bar{x}) = 0$, and $h_{22}(\bar{x}) = \lambda_2$. Moreover, we put 
\[\beta_1 = 1 - \psi'(\lambda_1-\lambda_2)\] 
and 
\[\beta_2 = 1 + \psi'(\lambda_1-\lambda_2).\] 
By assumption, we have $0 \leq \psi'(s) < 1$. This implies $0 < \beta_1 \leq \beta_2$.

\begin{lemma}
\label{reflection}
Let $\tau: \mathbb{R}^4 \to \mathbb{R}^4$ denote the reflection across the hyperplane orthogonal to $F(\bar{x})-F(\bar{y})$. Then $F(\bar{y}) = \tau(F(\bar{x}))$ and $\nu(\bar{y}) = \tau(\nu(\bar{x}))$. 
\end{lemma} 

\textbf{Proof.} 
The identity $F(\bar{y}) = \tau(F(\bar{x}))$ is trivial. In order to show that $\nu(\bar{y}) = \tau(\nu(\bar{x}))$, we choose a connected component $N$ of $S^3 \setminus F(\Sigma)$ with the property that $\nu$ is the outward-pointing unit normal vector field to $N$. Moreover, we define 
\[D = \{p \in S^3: \Phi(\bar{x}) \, (1 - \langle F(\bar{x}),p \rangle) + \langle \nu(\bar{x}),p \rangle < 0\}.\] 
Note that $D$ is a geodesic ball in $S^3$ with radius less than $\frac{\pi}{2}$. Moreover, it is easy to see that $\cos t \, F(\bar{x}) - \sin t \, \nu(\bar{x}) \in D$ if $t>0$ is sufficiently small. This shows that $D \cap N \neq \emptyset$. On the other hand, since the function $Z$ is nonnegative, we have $D \subset S^3 \setminus F(\Sigma)$. Putting these facts together, we conclude that $D \subset N$. 

Since $Z(\bar{x},\bar{y}) = 0$, the point $F(\bar{y})$ lies on the boundary of $D$ and on the boundary of $N$. Moreover, the outward-pointing unit normal vector to $N$ at the point $F(\bar{y})$ is $\nu(\bar{y})$, and the outward-pointing unit normal vector to $D$ at $F(\bar{y})$ is $\tau(F(\bar{x}))$. Since $D \subset N$, the two normal vectors must coincide. Consequently, $\nu(\bar{y}) = \tau(\nu(\bar{x}))$, as claimed. \\

Lemma \ref{reflection} directly implies that 
\[\text{\rm span} \Big \{ \frac{\partial F}{\partial y_1}(\bar{y}),\frac{\partial F}{\partial y_2}(\bar{y}) \Big \} = \text{\rm span} \Big \{ \tau(\frac{\partial F}{\partial x_1}(\bar{x})),\tau(\frac{\partial F}{\partial x_2}(\bar{x})) \Big \}.\] 
In particular, we can choose local coordinates $(y_1,y_2)$ around $\bar{y}$ so that 
\[\frac{\partial F}{\partial y_i}(\bar{y}) = \tau(\frac{\partial F}{\partial x_i}(\bar{x})).\] 

\begin{proposition} 
\label{second.derivatives.1}
We have 
\begin{align*} 
&\sum_{i=1}^2 \beta_i \, \frac{\partial^2 Z}{\partial x_i^2}(\bar{x},\bar{y}) \\ 
&= \bigg ( \sum_{i=1}^2 \beta_i \, \frac{\partial^2 \Phi}{\partial x_i^2}(\bar{x}) - 2 \sum_{i=1}^2 \frac{\beta_i}{\Phi(\bar{x})-\lambda_i} \, \Big ( \frac{\partial \Phi}{\partial x_i}(\bar{x}) \Big )^2 \\ 
&\hspace{8mm} + \sum_{i=1}^2 \beta_i \, (\lambda_i^2-1) \, \Phi(\bar{x}) - \sum_{i=1}^2 \beta_i \, \lambda_i \, (\Phi(\bar{x})^2 - 1) \bigg ) \cdot (1 - \langle F(\bar{x}),F(\bar{y}) \rangle) \\ 
&+ \sum_{i=1}^2 \beta_i \, (\Phi(\bar{x}) - \lambda_i). 
\end{align*} 
\end{proposition}

\textbf{Proof.} 
Using the Codazzi equations, we obtain 
\[\sum_{i=1}^2 \beta_i \, \frac{\partial}{\partial x_i} h_{ik}(\bar{x}) = \sum_{i=1}^2 \beta_i \, \frac{\partial}{\partial x_k} h_{ii}(\bar{x}) = \frac{\partial}{\partial x_k} \big ( \lambda_1+\lambda_2-\psi(\lambda_1-\lambda_2) \big ) = 0\] 
at the point $\bar{x}$. This implies 
\begin{align*} 
&\sum_{i=1}^2 \beta_i \, \frac{\partial^2 Z}{\partial x_i^2}(\bar{x},\bar{y}) \\ 
&= \sum_{i=1}^2 \beta_i \, \frac{\partial^2 \Phi}{\partial x_i^2}(\bar{x}) \, (1 - \langle F(\bar{x}),F(\bar{y}) \rangle) - 2 \sum_{i=1}^2 \beta_i \, \frac{\partial \Phi}{\partial x_i}(\bar{x}) \, \Big \langle \frac{\partial F}{\partial x_i}(\bar{x}),F(\bar{y}) \Big \rangle \\ 
&+ \sum_{i=1}^2 \beta_i \, \Phi(\bar{x}) \, \langle F(\bar{x}),F(\bar{y}) \rangle + \sum_{i=1}^2 \beta_i \, \lambda_i \, \Phi(\bar{x}) \, \langle \nu(\bar{x}),F(\bar{y}) \rangle \\ 
&- \sum_{i=1}^2 \beta_i \, \lambda_i \, \langle F(\bar{x}),F(\bar{y}) \rangle - \sum_{i=1}^2 \beta_i \, \lambda_i^2 \, \langle \nu(\bar{x}),F(\bar{y}) \rangle \\ 
&= \bigg ( \sum_{i=1}^2 \beta_i \, \frac{\partial^2 \Phi}{\partial x_i^2}(\bar{x}) + \sum_{i=1}^2 \beta_i \, \lambda_i^2 \, \Phi(\bar{x}) - \sum_{i=1}^2 \beta_i \, \lambda_i \, \Phi(\bar{x})^2 \bigg ) \, (1 - \langle F(\bar{x}),F(\bar{y}) \rangle) \\ 
&- 2 \sum_{i=1}^2 \beta_i \, \frac{\partial \Phi}{\partial x_i}(\bar{x}) \, \Big \langle \frac{\partial F}{\partial x_i}(\bar{x}),F(\bar{y}) \Big \rangle \\ 
&+ \sum_{i=1}^2 \beta_i \, \Phi(\bar{x}) \, \langle F(\bar{x}),F(\bar{y}) \rangle - \sum_{i=1}^2 \beta_i \, \lambda_i \, \langle F(\bar{x}),F(\bar{y}) \rangle 
\end{align*} 
at $\bar{x}$. This gives 
\begin{align*} 
&\sum_{i=1}^2 \beta_i \, \frac{\partial^2 Z}{\partial x_i^2}(\bar{x},\bar{y}) \\ 
&= \bigg ( \sum_{i=1}^2 \beta_i \, \frac{\partial^2 \Phi}{\partial x_i^2}(\bar{x}) + \sum_{i=1}^2 \beta_i \, (\lambda_i^2-1) \, \Phi(\bar{x}) - \sum_{i=1}^2 \beta_i \, \lambda_i \, (\Phi(\bar{x})^2-1) \bigg ) \\ 
&\hspace{8mm} \cdot (1 - \langle F(\bar{x}),F(\bar{y}) \rangle) \\ 
&- 2 \sum_{i=1}^2 \beta_i \, \frac{\partial \Phi}{\partial x_i}(\bar{x}) \, \Big \langle \frac{\partial F}{\partial x_i}(\bar{x}),F(\bar{y}) \Big \rangle + \sum_{i=1}^2 \beta_i \, \Phi(\bar{x}) - \sum_{i=1}^2 \beta_i \, \lambda_i.
\end{align*} 
Since 
\[\Big \langle \frac{\partial F}{\partial x_i}(\bar{x}),F(\bar{y}) \Big \rangle = \frac{1}{\Phi(\bar{x})-\lambda_i} \, \frac{\partial \Phi}{\partial x_i}(\bar{x}) \, (1-\langle F(\bar{x}),F(\bar{y}) \rangle),\] 
the assertion follows. \\

\begin{proposition}
\label{second.derivatives.2}
Assume that $\frac{\partial F}{\partial y_i}(\bar{y}) = \tau(\frac{\partial F}{\partial x_i}(\bar{x}))$. Then  
\[\frac{\partial^2 Z}{\partial x_i \, \partial y_i}(\bar{x},\bar{y}) = \lambda_i - \Phi(\bar{x}).\] 
In particular, 
\[\sum_{i=1}^2 \beta_i \, \frac{\partial^2 Z}{\partial x_i \, \partial y_i}(\bar{x},\bar{y}) = \sum_{i=1}^2 \beta_i \, (\lambda_i - \Phi(\bar{x})).\] 
\end{proposition}

\textbf{Proof.} 
The proof of Proposition 5 in \cite{Brendle1} goes through unchanged. \\

\begin{proposition} 
\label{second.derivatives.3}
We have 
\[\sum_{i=1}^2 \beta_i \, \frac{\partial^2 Z}{\partial y_i^2}(\bar{x},\bar{y}) \leq \sum_{i=1}^2 \beta_i \, (\Phi(\bar{x}) - \lambda_i).\] 
\end{proposition}

\textbf{Proof.} 
We have 
\begin{align*} 
\sum_{i=1}^2 \beta_i \, \frac{\partial^2 Z}{\partial y_i^2}(\bar{x},\bar{y}) 
&= \sum_{i=1}^2 \beta_i \, \Big \langle \nu(\bar{x})-\Phi(\bar{x}) \, F(\bar{x}),\frac{\partial^2 F}{\partial y_i^2}(\bar{y}) \Big \rangle \\ 
&= -\sum_{i=1}^2 \beta_i \, \Big \langle \nu(\bar{x})-\Phi(\bar{x}) \, F(\bar{x}),F(\bar{y}) \Big \rangle \\ 
&- \sum_{i=1}^2 \beta_i \, h_{ii}(\bar{y}) \, \Big \langle \nu(\bar{x})-\Phi(\bar{x}) \, F(\bar{x}),\nu(\bar{y}) \Big \rangle. 
\end{align*} 
By Lemma \ref{reflection}, we have 
\[\nu(\bar{x}) - \Phi(\bar{x}) \, (F(\bar{x})-F(\bar{y})) = \tau(\nu(\bar{x})) = \nu(\bar{y}),\] 
hence 
\[\sum_{i=1}^2 \beta_i \, \frac{\partial^2 Z}{\partial y_i^2}(\bar{x},\bar{y}) = \sum_{i=1}^2 \beta_i \, (\Phi(\bar{x}) - h_{ii}(\bar{y})).\] 
Let $\mu_1 \geq \mu_2$ denote the principal curvatures at the point $\bar{y}$. Since $\psi$ is convex, we have 
\begin{align*} 
(\mu_1-\lambda_1)+(\mu_2-\lambda_2) 
&= \psi(\mu_1-\mu_2) - \psi(\lambda_1-\lambda_2) \\ 
&\geq \psi'(\lambda_1-\lambda_2) \, ((\mu_1-\lambda_1)-(\mu_2-\lambda_2)), 
\end{align*}
hence 
\[\sum_{i=1}^2 \beta_i \, (\mu_i-\lambda_i) \geq 0.\] 
Since $\beta_1 \leq \beta_2$, we conclude that 
\[\sum_{i=1}^2 \beta_i \, (h_{ii}(\bar{y}) - \lambda_i) \geq \sum_{i=1}^2 \beta_i \, (\mu_i-\lambda_i) \geq 0.\] 
Putting these facts together, we conclude that 
\[\sum_{i=1}^2 \beta_i \, \frac{\partial^2 Z}{\partial y_i^2}(\bar{x},\bar{y}) \leq \sum_{i=1}^2 \beta_i \, (\Phi(\bar{x}) - \lambda_i),\] 
as claimed. \\

Combining Proposition \ref{second.derivatives.1}, Proposition \ref{second.derivatives.2}, and Proposition \ref{second.derivatives.3}, we can draw the following conclusion:

\begin{corollary} 
\label{key.ingredient}
We have
\begin{align*} 
0 &\leq \sum_{i=1}^2 \beta_i \, \frac{\partial^2 Z}{\partial x_i^2}(\bar{x},\bar{y}) + 2 \sum_{i=1}^2 \beta_i \, \frac{\partial^2 Z}{\partial x_i \, \partial y_i}(\bar{x},\bar{y}) + \sum_{i=1}^2 \beta_i \, \frac{\partial^2 Z}{\partial y_i^2}(\bar{x},\bar{y}) \\ 
&\leq \bigg ( \sum_{i=1}^2 \beta_i \, \frac{\partial^2 \Phi}{\partial x_i^2}(\bar{x}) - 2 \sum_{i=1}^2 \frac{\beta_i}{\Phi(\bar{x})-\lambda_i} \, \Big ( \frac{\partial \Phi}{\partial x_i}(\bar{x}) \Big )^2 \\ 
&\hspace{8mm} + \sum_{i=1}^2 \beta_i \, (\lambda_i^2-1) \, \Phi(\bar{x}) - \sum_{i=1}^2 \beta_i \, \lambda_i \, (\Phi(\bar{x})^2 - 1) \bigg ) \, (1 - \langle F(\bar{x}),F(\bar{y}) \rangle). 
\end{align*} 
\end{corollary}

\section{The choice of the function $\Phi$}

In this section, we construct a function $\Phi$ such that 
\[\sum_{i=1}^2 \beta_i \, D_{i,i}^2 \Phi - 2 \sum_{i=1}^2 \frac{\beta_i}{\Phi-\lambda_i} \, (D_i \Phi)^2 + \sum_{i=1}^2 \beta_i \, (\lambda_i^2-1) \, \Phi - \sum_{i=1}^2 \beta_i \, \lambda_i \, (\Phi^2 - 1) < 0,\] 
where the coefficients $\beta_i$ are defined as above. To that end, we need the following result, which can be viewed as an analogue of the Simons identity (cf. \cite{Simons}) for the norm of the second fundamental form of a minimal surface: 

\begin{proposition} 
\label{pde.for.eigenvalues}
Let $\lambda_1$ and $\lambda_2$ be the principal curvatures (viewed as functions on $\Sigma$). Then 
\begin{align*} 
&\sum_{i=1}^2 \beta_i \, D_{i,i}^2 \lambda_1 + \sum_{i=1}^2 \beta_i \, (\lambda_i^2-1) \, \lambda_1 - \sum_{i=1}^2 \beta_i \, \lambda_i \, (\lambda_1^2-1) \\ 
&= 2 \, \frac{\beta_2}{\lambda_1-\lambda_2} \, ((D_1 \lambda_2)^2 + (D_2 \lambda_1)^2) + \psi''(\lambda_1-\lambda_2) \, \frac{(\beta_1+\beta_2)^2}{\beta_1^2} \, (D_1 \lambda_2)^2 
\end{align*} 
and 
\begin{align*} 
&\sum_{i=1}^2 \beta_i \, D_{i,i}^2 \lambda_2 + \sum_{i=1}^2 \beta_i \, (\lambda_i^2-1) \, \lambda_2 - \sum_{i=1}^2 \beta_i \, \lambda_i \, (\lambda_2^2-1) \\ 
&= -2 \, \frac{\beta_1}{\lambda_1-\lambda_2} \, ((D_1 \lambda_2)^2 + (D_2 \lambda_1)^2) + \psi''(\lambda_1-\lambda_2) \, \frac{(\beta_1+\beta_2)^2}{\beta_2^2} \, (D_2 \lambda_1)^2. 
\end{align*} 
\end{proposition}

\textbf{Proof.} 
We first relate the Hessian of the function $\lambda_i$ to the second covariant derivatives of the second fundamental form. The Hessian of the function $\lambda_i$ is given by 
\[D_{k,k}^2 \lambda_i = D_{k,k}^2 h_{ii} - (-1)^i \, \frac{2}{\lambda_1 - \lambda_2} \, (D_k h_{12})^2.\] 
Moreover, we have the commutator identity 
\[D_{k,k}^2 h_{ii} - D_{i,i}^2 h_{kk} = D_{k,i}^2 h_{ik} - D_{i,k}^2 h_{ik} = (\lambda_i^2-1) \, \lambda_k - \lambda_i \, (\lambda_k^2-1).\] 
Differentiating the relation $\lambda_1+\lambda_2 - \psi(\lambda_1-\lambda_2) = 0$ gives 
\begin{align*} 
0 &= D_{k,k}^2 (\lambda_1+\lambda_2 - \psi(\lambda_1-\lambda_2)) \\ 
&= \sum_{i=1}^2 \beta_i \, D_{k,k}^2 \lambda_i - \psi''(\lambda_1-\lambda_2) \, (D_k \lambda_1- D_k \lambda_2)^2 \\ 
&= \sum_{i=1}^2 \beta_i \, D_{k,k}^2 h_{ii} + 2 \, \frac{\beta_1-\beta_2}{\lambda_1-\lambda_2} \, (D_k h_{12})^2 - \psi''(\lambda_1-\lambda_2) \, (D_k \lambda_1- D_k \lambda_2)^2 \\ 
&= \sum_{i=1}^2 \beta_i \, D_{i,i}^2 h_{kk} + \sum_{i=1}^2 \beta_i \, (\lambda_i^2-1) \, \lambda_k - \sum_{i=1}^2 \beta_i \, \lambda_i \, (\lambda_k^2 - 1) \\ 
&+ 2 \, \frac{\beta_1-\beta_2}{\lambda_1-\lambda_2} \, (D_k h_{12})^2 - \psi''(\lambda_1-\lambda_2) \, (D_k \lambda_1- D_k \lambda_2)^2 \\ 
&= \sum_{i=1}^2 \beta_i \, D_{i,i}^2 \lambda_k + \sum_{i=1}^2 \beta_i \, (\lambda_i^2-1) \, \lambda_k - \sum_{i=1}^2 \beta_i \, \lambda_i \, (\lambda_k^2 - 1) \\ 
&+ 2 \, (-1)^k \sum_{i=1}^2 \frac{\beta_i}{\lambda_1-\lambda_2} \, (D_i h_{12})^2 + 2 \, \frac{\beta_1-\beta_2}{\lambda_1-\lambda_2} \, (D_k h_{12})^2 \\ 
&- \psi''(\lambda_1-\lambda_2) \, (D_k \lambda_1- D_k \lambda_2)^2. 
\end{align*} 
This gives 
\begin{align*} 
&\sum_{i=1}^2 \beta_i \, D_{i,i}^2 \lambda_1 + \sum_{i=1}^2 \beta_i \, (\lambda_i^2-1) \, \lambda_1 - \sum_{i=1}^2 \beta_i \, \lambda_i \, (\lambda_1^2-1) \\ 
&= 2 \, \frac{\beta_2}{\lambda_1-\lambda_2} \, ((D_1 \lambda_2)^2 + (D_2 \lambda_1)^2) + \psi''(\lambda_1-\lambda_2) \, (D_1 \lambda_1 - D_1 \lambda_2)^2 
\end{align*} 
and 
\begin{align*} 
&\sum_{i=1}^2 \beta_i \, D_{i,i}^2 \lambda_2 + \sum_{i=1}^2 \beta_i \, (\lambda_i^2-1) \, \lambda_2 - \sum_{i=1}^2 \beta_i \, \lambda_i \, (\lambda_2^2-1) \\ 
&= -2 \, \frac{\beta_1}{\lambda_1-\lambda_2} \, ((D_1 \lambda_2)^2 + (D_2 \lambda_1)^2) + \psi''(\lambda_1-\lambda_2) \, (D_2 \lambda_1 - D_2 \lambda_2)^2. 
\end{align*} 
The assertion now follows from the fact that $\sum_{i=1}^2 \beta_i \, D_k \lambda_i = 0$ for $k=1,2$. \\

\begin{proposition} 
\label{barrier}
Let $\Phi = \alpha\lambda_1 - (\alpha-1) \lambda_2$, where $\lambda_1 \geq \lambda_2$ are the principal curvatures and $\alpha > 1$ is a constant. Then 
\[\sum_{i=1}^2 \beta_i \, D_{i,i}^2 \Phi - 2 \sum_{i=1}^2 \frac{\beta_i}{\Phi-\lambda_i} \, (D_i \Phi)^2 + \sum_{i=1}^2 \beta_i \, (\lambda_i^2-1) \, \Phi - \sum_{i=1}^2 \beta_i \, \lambda_i \, (\Phi^2 - 1) < 0.\] 
\end{proposition}

\textbf{Proof.}
It follows from Proposition \ref{pde.for.eigenvalues} that 
\begin{align*} 
&\sum_{i=1}^2 \beta_i \, D_{i,i}^2 \Phi + \sum_{i=1}^2 \beta_i \, (\lambda_i^2-1) \, \Phi - \sum_{i=1}^2 \beta_i \, \lambda_i \, (\Phi^2 - 1) \\ 
&= -\alpha(\alpha-1) \sum_{i=1}^2 \beta_i \, \lambda_i \, (\lambda_1-\lambda_2)^2 + 2 \, \frac{(\alpha-1)\beta_1+\alpha\beta_2}{\lambda_1-\lambda_2} \, ((D_1\lambda_2)^2+(D_2\lambda_1)^2) \\ 
&+ \alpha \, \psi''(\lambda_1-\lambda_2) \, \frac{(\beta_1+\beta_2)^2}{\beta_1^2} \, (D_1 \lambda_2)^2 - (\alpha-1) \, \psi''(\lambda_1-\lambda_2) \, \frac{(\beta_1+\beta_2)^2}{\beta_2^2} \, (D_2 \lambda_1)^2. 
\end{align*}
Moreover, using the relation $\sum_{i=1}^2 \beta_i \, D_k \lambda_i = 0$, we obtain 
\[D_1 \Phi = \alpha \, D_1 \lambda_1 - (\alpha-1) \, D_1 \lambda_2 = -\frac{(\alpha-1) \beta_1 + \alpha\beta_2}{\beta_1} \, D_1 \lambda_2\] 
and 
\[D_2 \Phi = \alpha \, D_2 \lambda_1 - (\alpha-1) \, D_2 \lambda_2 = \frac{(\alpha-1) \beta_1 + \alpha \beta_2}{\beta_2} \, D_2 \lambda_1.\] 
This implies 
\begin{align*} 
&\sum_{i=1}^2 \frac{\beta_i}{\Phi-\lambda_i} \, (D_i \Phi)^2 \\ 
&= \frac{1}{\alpha-1} \, \frac{\beta_1}{\lambda_1-\lambda_2} \, (D_1 \Phi)^2 + \frac{1}{\alpha} \, \frac{\beta_2}{\lambda_1-\lambda_2} \, (D_2 \Phi)^2 \\ 
&= \frac{(\alpha-1) \beta_1 + \alpha\beta_2}{\lambda_1-\lambda_2} \, \Big ( \frac{(\alpha-1) \beta_1 + \alpha\beta_2}{(\alpha-1) \beta_1} \, (D_1 \lambda_2)^2 + \frac{(\alpha-1) \beta_1 + \alpha\beta_2}{\alpha \beta_2} \, (D_2 \lambda_1)^2 \Big ). 
\end{align*} 
Putting these facts together gives 
\begin{align*} 
&\sum_{i=1}^2 \beta_i \, D_{i,i}^2 \Phi - 2 \sum_{i=1}^2 \frac{\beta_i}{\Phi-\lambda_i} \, (D_i \Phi)^2 + \sum_{i=1}^2 \beta_i \, (\lambda_i^2-1) \, \Phi - \sum_{i=1}^2 \beta_i \, \lambda_i \, (\Phi^2 - 1) \\ 
&= -\alpha(\alpha-1) \sum_{i=1}^2 \beta_i \, \lambda_i \, (\lambda_1-\lambda_2)^2 \\ 
&- 2 \, \frac{(\alpha-1) \beta_1 + \alpha\beta_2}{\lambda_1-\lambda_2} \, \Big ( \frac{\alpha\beta_2}{(\alpha-1) \beta_1} \, (D_1 \lambda_2)^2 + \frac{(\alpha-1) \beta_1}{\alpha \beta_2} \, (D_2 \lambda_1)^2 \Big ) \\ 
&+ \alpha \, \psi''(\lambda_1-\lambda_2) \, \frac{(\beta_1+\beta_2)^2}{\beta_1^2} \, (D_1 \lambda_2)^2 - (\alpha-1) \, \psi''(\lambda_1-\lambda_2) \, \frac{(\beta_1+\beta_2)^2}{\beta_2^2} \, (D_2 \lambda_1)^2. 
\end{align*} 
Since $\psi''(s) \geq 0$, the last term on the right hand side has a favorable sign. To estimate the second to last term, we use the inequality $\psi''(s) \leq \frac{1-\psi'(s)^2}{s}$. This gives 
\[\psi''(\lambda_1-\lambda_2) \leq \frac{1-\psi'(\lambda_1-\lambda_2)^2}{\lambda_1-\lambda_2} = \frac{\beta_1\beta_2}{\lambda_1-\lambda_2} \, \frac{2}{\beta_1+\beta_2}.\] 
Thus, we conclude that 
\begin{align*} 
&\sum_{i=1}^2 \beta_i \, D_{i,i}^2 \Phi - 2 \sum_{i=1}^2 \frac{\beta_i}{\Phi-\lambda_i} \, (D_i \Phi)^2 + \sum_{i=1}^2 \beta_i \, (\lambda_i^2-1) \, \Phi - \sum_{i=1}^2 \beta_i \, \lambda_i \, (\Phi^2 - 1) \\ 
&\leq -\alpha(\alpha-1) \sum_{i=1}^2 \beta_i \, \lambda_i \, (\lambda_1-\lambda_2)^2 \\ 
&- 2 \, \frac{\beta_2}{\lambda_1-\lambda_2} \, \frac{\alpha\beta_2}{(\alpha-1) \beta_1} \, (D_1 \lambda_2)^2 - 2 \, \frac{(\alpha-1) \beta_1 + \alpha\beta_2}{\lambda_1-\lambda_2} \, \frac{(\alpha-1) \beta_1}{\alpha \beta_2} \, (D_2 \lambda_1)^2 \\ 
&< 0.
\end{align*} 
In the last step we have used the inequality $\sum_{i=1}^2 \beta_i \, \lambda_i > 0$. This follows from our assumption that $s \, \psi'(s) < \psi(s)$. \\

\section{Proof of Theorem \ref{main.result}}

In this section, we complete the proof of Theorem \ref{main.result}. As above, we assume that $F: \Sigma \to S^3$ is an embedding which satisfies the equation $\lambda_1+\lambda_2 = \psi(\lambda_1-\lambda_2)$, where $\psi$ satisfies the structure conditions from Theorem \ref{main.result}. 

\begin{proposition}
\label{inscribed.radius}
We have 
\[\lambda_1(x) \, (1 - \langle F(x),F(y) \rangle) + \langle \nu(x),F(y) \rangle \geq 0\] 
for all $x,y \in \Sigma$, where $\lambda_1(x)$ denotes the larger one of the principal curvatures at $x$. In other words, at each point on $\Sigma$ the inscribed radius equals the curvature radius.
\end{proposition}

\textbf{Proof.} 
For each $\alpha>1$, we define  
\[Z_\alpha(x,y) = \Phi_\alpha(x) \, (1 - \langle F(x),F(y) \rangle) + \langle \nu(x),F(y) \rangle,\] 
where $\Phi_\alpha(x) = \alpha \, \lambda_1(x) - (\alpha-1) \, \lambda_2(x)$. Since $F$ has no umbilical points, the function $Z_\alpha$ is nonnegative when $\alpha$ is sufficiently large. Let now 
\[\kappa = \inf \{\alpha>1: \text{\rm $Z_\alpha(x,y) \geq 0$ for all $x,y \in \Sigma$}\}.\] 
Clearly, $Z_\kappa(x,y) \geq 0$ for all points $x,y \in \Sigma$. 

We claim that $\kappa=1$. Indeed, if $\kappa > 1$, then we can find a pair of points $\bar{x},\bar{y} \in \Sigma$ such that $\bar{x} \neq \bar{y}$ and $Z_\kappa(\bar{x},\bar{y}) = 0$. By Corollary \ref{key.ingredient}, we have 
\begin{align*} 
0 &\leq \sum_{i=1}^2 \beta_i \, \frac{\partial^2 Z_\kappa}{\partial x_i^2}(\bar{x},\bar{y}) + 2 \sum_{i=1}^2 \beta_i \, \frac{\partial^2 Z_\kappa}{\partial x_i \, \partial y_i}(\bar{x},\bar{y}) + \sum_{i=1}^2 \beta_i \, \frac{\partial^2 Z_\kappa}{\partial y_i^2}(\bar{x},\bar{y}) \\ 
&\leq \bigg ( \sum_{i=1}^2 \beta_i \, \frac{\partial^2 \Phi_\kappa}{\partial x_i^2}(\bar{x}) - 2 \sum_{i=1}^2 \frac{\beta_i}{\Phi_\kappa(\bar{x})-\lambda_i} \, \Big ( \frac{\partial \Phi_\kappa}{\partial x_i}(\bar{x}) \Big )^2 \\ 
&\hspace{8mm} + \sum_{i=1}^2 \beta_i \, (\lambda_i^2-1) \, \Phi_\kappa(\bar{x}) - \sum_{i=1}^2 \beta_i \, \lambda_i \, (\Phi_\kappa(\bar{x})^2 - 1) \bigg ) \, (1 - \langle F(\bar{x}),F(\bar{y}) \rangle). 
\end{align*} 
This contradicts Proposition \ref{barrier}. Thus, we conclude that $\kappa=1$. This implies that 
\[Z_1(x,y) = \lambda_1(x) \, (1 - \langle F(x),F(y) \rangle) + \langle \nu(x),F(y) \rangle \geq 0\] 
for all points $x,y \in \Sigma$. This completes the proof of Proposition \ref{inscribed.radius}. \\

Using Proposition \ref{inscribed.radius}, we are able to show that the principal curvatures are constant along one family of curvature lines:

\begin{corollary}
We have $D_1 \lambda_1 = D_1 \lambda_2 = 0$ at each point on $\Sigma$.
\end{corollary}

\textbf{Proof.} 
Let us fix a point $\bar{x} \in \Sigma$, and let $\{e_1,e_2\}$ be an orthonormal basis of $T_{\bar{x}} \Sigma$ such that $h(e_1,e_1) = \lambda_1(\bar{x})$, $h(e_1,e_2) = 0$, and $h(e_2,e_2) = \lambda_2(\bar{x})$. For abbreviation, we put 
\[\xi = \lambda_1(\bar{x}) \, F(\bar{x}) - \nu(\bar{x}) \in \mathbb{R}^4.\] 
Let us consider a geodesic $\sigma: \mathbb{R} \to \Sigma$ be a geodesic on $\Sigma$ satisfying $\sigma(0) = \bar{x}$ and $\sigma'(0) = e_1$. By Proposition \ref{inscribed.radius}, the function 
\[f(t) = \lambda_1(\bar{x}) - \langle \xi,F(\sigma(t)) \rangle\] 
is nonnegative for all $t$. Moreover, the derivatives of $f(t)$ are given by  
\[f'(t) = -\langle \xi,dF_{\sigma(t)}(\sigma'(t)) \rangle,\]
\[f''(t) = \langle \xi,F(\sigma(t)) \rangle + h(\sigma'(t),\sigma'(t)) \, \langle \xi,\nu(\sigma(t)) \rangle,\] 
and 
\begin{align*} 
f'''(t) 
&= \langle \xi,dF_{\sigma(t)}(\sigma'(t)) \rangle + h(\sigma'(t),\sigma'(t)) \, \langle \xi,D_{\sigma'(t)} \nu \rangle \\ 
&+ (D_{\sigma'(t)}^\Sigma h)(\sigma'(t),\sigma'(t)) \, \langle \xi,\nu(\sigma(t)) \rangle. 
\end{align*} 
Putting $t=0$, we obtain $f(0) = f'(0) = f''(0) = 0$. Since the function $f(t)$ is nonnegative, we conclude that $f'''(0) = 0$. This gives $(D_{e_1}^\Sigma h)(e_1,e_1) = 0$. Consequently, we have $D_1 \lambda_1 = 0$ at the point $\bar{x}$. Since $\sum_{i=1}^2 \beta_i \, D_k \lambda_i = 0$ for $k=1,2$, it follows that $D_1 \lambda_2 = 0$ at the point $\bar{x}$. \\

We now complete the proof of Theorem \ref{main.result}. Let us choose a non-vanishing function $\varphi$  such that 
\[\varphi'(s) = -\frac{1+\psi'(s)}{2s} \, \varphi(s)\] 
for $s>0$. Moreover, we define a vector field $V$ on $\Sigma$ by 
\[V = \varphi(\lambda_1-\lambda_2) \, e_1.\]

\begin{lemma}
\label{lie.bracket}
The vector field $V$ satisfies $[V,e_1] = [V,e_2] = 0$.
\end{lemma}

\textbf{Proof.} 
Using the identity $D_1 (\lambda_1-\lambda_2) = 0$, we obtain $[V,e_1] = 0$. Moreover, the identity $D_1 \lambda_2 = 0$ implies $D_{e_2}^\Sigma e_1 = 0$. Since 
\[D_{e_1}^\Sigma e_2 = -\frac{1}{\lambda_1-\lambda_2} \, D_2 \lambda_1 \, e_1,\] 
we conclude that 
\[[e_1,e_2] = D_{e_1}^\Sigma e_2 - D_{e_2}^\Sigma e_1 = -\frac{1}{\lambda_1-\lambda_2} \, D_2 \lambda_1 \, e_1.\] 
Consequently, we have 
\begin{align*} 
[V,e_2] 
&= \varphi(\lambda_1-\lambda_2) \, [e_1,e_2] - \varphi'(\lambda_1-\lambda_2) \, D_2 (\lambda_1-\lambda_2) \, e_1 \\ 
&= -\frac{1}{\lambda_1-\lambda_2} \, \varphi(\lambda_1-\lambda_2) \, D_2 \lambda_1 \, e_1 - \varphi'(\lambda_1-\lambda_2) \, \frac{2}{1+\psi'(\lambda_1-\lambda_2)} \, D_2 \lambda_1 \, e_1 \\ 
&= 0. 
\end{align*}
This completes the proof of Lemma \ref{lie.bracket}. \\

Since $[V,e_1] = [V,e_2] = 0$ and $V(\lambda_1) = V(\lambda_2) = 0$, we obtain 
\[(\mathscr{L}_V g)(e_i,e_j) = V(g(e_i,e_j)) = 0\] 
and 
\[(\mathscr{L}_V h)(e_i,e_j) = V(h(e_i,e_j)) = 0.\] 
Thus, $\mathscr{L}_V g = \mathscr{L}_V h = 0$. Using the normal exponential map, we may extend $V$ to a Killing vector field which is defined on a small tubular neighborhood of $F(\Sigma)$. Consequently, the vector field $\xi$ must be the restriction of an ambient Killing vector field on $S^3$. In other words, there exists an anti-symmetric matrix $Q \in \mathfrak{so}(4)$ such that $V(x) = Q \, F(x)$ for all $x \in \Sigma$. 

It remains to show that the matrix $Q$ has rank $2$. To see this, we differentiate the relation $V(x) = Q \, F(x)$ along $e_2$. Since $g(V,e_2) = h(V,e_2) = 0$, we obtain $D_{e_2}^\Sigma V = Q \, e_2$. On the other hand, we have $D_{e_2}^\Sigma V = -\frac{\varphi(\lambda_1-\lambda_2)}{\lambda_1-\lambda_2} \, D_2 \lambda_1 \, e_1$. Thus, 
\[Q \, \Big ( e_2 + \frac{1}{\lambda_1-\lambda_2} \, D_2 \lambda_1 \, F(x) \Big ) = Q \, e_2 + \frac{\varphi(\lambda_1-\lambda_2)}{\lambda_1-\lambda_2} \, D_2 \lambda_1 \, e_1 = 0.\] 
This shows that $Q$ has non-trivial nullspace. Consequently, $Q$ has rank $2$, as claimed.

\end{document}